\documentclass[11pt,leqno]{article}

\usepackage{amsmath,amsfonts,amscd,amssymb,theorem}

\long\def\comment#1\endcomment{}

\makeatletter
\begingroup
\gdef\th@dotted{\normalfont\itshape
  \def\@begintheorem##1##2{%
        \item[\hskip\labelsep \theorem@headerfont ##1\ ##2.]}%
\def\@opargbegintheorem##1##2##3{%
   \item[\hskip\labelsep \theorem@headerfont ##1\ ##2\ (##3).]}}
\endgroup
\makeatother

\theoremstyle{dotted}

\newtheorem{theorem}{Theorem}[section]
\newtheorem{lemma}[theorem]{Lemma}

\newtheorem{prop}[theorem]{Proposition}

\makeatletter
\begingroup
\gdef\th@upshape{\normalfont
  \def\@begintheorem##1##2{%
        \item[\hskip\labelsep \theorem@headerfont ##1\ ##2.]}%
\def\@opargbegintheorem##1##2##3{%
   \item[\hskip\labelsep \theorem@headerfont ##1\ ##2\ (##3).]}}
\endgroup
\makeatother

\theoremstyle{upshape}

\newtheorem{exa}[theorem]{Example}

\makeatletter
\renewcommand{\subsection}{\@startsection{subsection}{2}{0pt}{-3ex
plus -1ex minus -0.2ex}{-2mm plus -0pt minus
-2pt}{\normalfont\bfseries}} 
\renewcommand{\subsubsection}{\@startsection{subsubsection}{3}{0pt}{-3ex
plus -1ex minus -0.2ex}{-2mm plus -0pt minus
-2pt}{\normalfont\bfseries}} 
\makeatother

\makeatletter
\@addtoreset{equation}{section}
\makeatother

\newcommand{\cntrct}                
{\hspace{2pt}\raisebox{1pt}{\text{$\lrcorner$}}\hspace{2pt}}

\newcommand{\proof}[1][Proof.]{\smallskip\noindent{\em #1}}
\def\endproof{\hfill\ensuremath{\square}\par\medskip}

\renewcommand{\labelenumi}{{\normalfont(\roman{enumi})}}

\def\eqref#1{\thetag{\ref{#1}}}

\let\latexref=\ref
\def\ref#1{{\normalfont{\latexref{#1}}}}

\newcommand{\idot}{{\:\raisebox{1pt}{\text{\circle*{1.5}}}}}
\newcommand{\hdot}{{\:\raisebox{3pt}{\text{\circle*{1.5}}}}}
\newcommand{\eps}{\varepsilon}
\renewcommand{\phi}{\varphi}

\newcommand{\vH}{\check{H}}

\newcommand{\Ker}{\operatorname{Ker}}

\newcommand{\id}{\operatorname{\sf id}}
\newcommand{\Id}{\operatorname{\sf Id}}

\newcommand{\tr}{\operatorname{\sf tr}}
\newcommand{\colim}{\operatorname{colim}}
\renewcommand{\lim}{\operatorname{lim}}

\newcommand{\End}{{\operatorname{End}}}

\newcommand{\amod}{{\text{\rm -mod}}}

\newcommand{\ppt}{{\sf pt}}

\newcommand{\cchar}{\operatorname{\sf char}}

\newcommand{\Z}{{\mathbb Z}}

\newcommand{\Fun}{\operatorname{Fun}}

\newcommand{\E}{\mathcal{E}}

\newcommand{\Stab}{\operatorname{Stab}}

\newcommand{\Top}{\operatorname{Top}}
\newcommand{\Ho}{\operatorname{Ho}}

\newcommand{\copr}{\sqcup}

\newcommand{\Tr}{\operatorname{Tr}}

\newcommand{\ev}{\operatorname{\sf ev}}

\newcommand{\F}{\operatorname{\mathbb{F}}}

\newcommand{\Pp}{\operatorname{\mathbb{P}}}

\newcommand{\I}{\operatorname{\sf I}}

\newcommand{\ssl}{\mathfrak{s}\mathfrak{l}}
\newcommand{\bfr}{\mathfrak{b}}
\newcommand{\ufr}{\mathfrak{u}}
\newcommand{\6}{\partial}
\newcommand{\calo}{\mathcal{O}}

\newcommand{\gfr}{\mathfrak{g}}

\title{B\"okstedt periodicity generator via $K$-theory}

\author{A. Fonarev and D. Kaledin\thanks{Both authors supported by
    the Basis Foundation, grant 18-1-6-95-1, Leader (Math).}}

\begin{document}

\maketitle

\tableofcontents

\section*{Introduction.}

In its most basic form, ``B\"okstedt periodicity'' refers to the
computation of the Topological Hochschild Homology $THH_\idot(k)$ of
a prime field $k=\F_p$ with $p$ elements; the answer is that
$THH_\idot(k) \cong k[v]$ is the polynomial algebra in one generator
$v$ of homological degree $2$. The computation was done by
B\"okstedt \cite{bo} in the same paper where Topological Hochschild
Homology was first defined. Slightly later, in a wonderful series of
papers \cite{P.et.al}, \cite{JP}, \cite{PW} stemming from
\cite{P.0}, Pirashvili and others showed that $THH_\idot(k)$ can be
identified with several other homology theories including Mac Lane
Homology $HM_\idot(k)$, where the same periodicity has been earlier
established by Breen \cite{breen}. A good overview of the whole
story as it was seen in 1998 can be found in \cite{LP}, and a more
recent survey is available in \cite{HN}.

As of now, B\"okstedt periodicity has numerous applications,
including some very prominent ones such as \cite{BMS}, and several
proofs none of which are elementary. The latter is somewhat
surprising, given that the question is so fundamental, and the
answer is so miraculously simple. What happens is, one can easily
construct a spectral sequence converging to $THH(k)$ whose
$E_2$-term can be read off the dual Steenrod algebra, and the
sequence degenerates at $E_p$. However, if $p$ is not $2$, then the
$E_2$-term is rather large, and while it already contains $v$, there
are many extra elements that are killed off by the differential
$d_{p-1}$. So, the spectral sequence is highly non-degenerate. One
has to prove this somehow, and then show that $THH(k)\cong k[v]$ as
a $k$-algebra (this is also non-trivial since in the $E_p$ term,
$v^p=0$). Roughly speaking, the existing proofs belong to one of the
following two types.
\begin{enumerate}
\item ``Topological'' proofs that use Dyer-Lashof operations. This
  was the original approach of B\"okstedt \cite{bo}. Recently, a
  pretty definitive proof along the same lines was given in
  \cite{KN} that reduces the claim to the following amazing one-line
  observation: the dual Steenrod algebra for $k$ is freely generated
  by one element of degree $1$ as an $E_2$-algebra. This yields
  B\"okstedt periodicity almost immediately, and the observation
  itself follows from classical computations of Dyer-Lashof
  operations.
\item ``Representation-theoretic'' proofs. These use the
  identification between $THH$ and Mac Lane Homology, and then work
  with Mac Lane Cohomology instead of Homology, and use polynomial
  functors and modular representation theory of $GL(n)$. The
  original proof along these lines is \cite{slf}, and it can be
  considerably simplified using more recent technology (see
  e.g. \cite[Section 11]{K.bb}).
\end{enumerate}
In our opinion, none of these proofs is so satisfactory as to make
the other ones obsolete. If one just wants shortness, it is hard to
beat \cite{KN}. However, while the argument proves the theorem, it
does not explain the miracle but moves it one step up. Why would the
dual Steenrod algebra, an $E_\infty$-al\-ge\-bra by nature, become
so simple when considered as an $E_n$-algebra for $n=2$? One has a
feeling that there should be a categorical reason for this, but at
the moment, it is not clear, at least to us, what it might be. Also,
the duality between Mac Lane Homology and Mac Lane Cohomology is not
visible in the topological approach, and as far as we know, a
possible cohomological counterpart to $THH$ has not been properly
studied.

In the representation-theoretic approach, it is the cohomology that
plays the main role, and the non-degeneration of the spectral
sequence acquires some sort of an explanation. At the end of the
day, the main fact is that the category of strictly polynomial
functors of some fixed degree has finite homological dimension. This
line of reasoning actually leads quite far --- for example, to the
celebrated finiteness results of Friedlander and Suslin \cite{FS}
and more recent strong results of Touz\'e \cite{tou}. However, there
is a price to pay: while non-degeneration for cohomology is
equivalent to non-degeneration for homology, this approach says
absolutely nothing about multiplication in homology. We obtain a
graded vector space isomorphism $THH_\idot(k) \cong k[v]$, but to
check that it is multiplicative, one needs an additional argument.

One possible argument goes as follows. Firstly, for dimension
reasons, it suffices to check that the generator $v$ in $THH_2(k)$
is not nilpotent. To do this, it suffices to construct a
multiplicative map $\phi_\idot:THH_\idot(k) \to R_\idot(k)$ to some
graded algebra $R_\idot(k)$ such that $\phi_\idot(v)^n \neq 0$ for
all $n \geq 0$. A tempting option is to take as $R_\idot(k)$ the
periodic cyclic homology $HP_\idot(k)$, with the multiplicative map
$\phi_\idot:THH_\idot(k) \to HP_\idot(k)$ extracted from the
cyclotomic structure map on $THH_\idot(k)$ of \cite{NS} by the
method suggested in \cite{mathew}. We have $HP_\idot(k) \cong
k((u))$, for a generator $u$ of cohomological degree $2$, and {\em a
  posteriori}, $\phi(v) = u^{-1}$; were we to know this {\em a
  priori}, this would be it.

The argument in \cite[Section 11]{K.bb} goes along the same lines
but in a slightly different way. Instead of $HP_\idot(k)$, one takes
$R_\idot(k) = \vH^\hdot(\Z/p\Z,k)$, the Tate cohomology algebra of
the cyclic group $\Z/p\Z$ with coefficients in the trivial module
$k$. It also has a generator $u$ of cohomological degree $2$. To
construct a multiplicative map $\phi_\idot$, one uses yet another
re-intrepretation of $THH_\idot(k)$ --- namely, the identification,
essentially due to \cite{JP}, between $THH_\idot(A,M)$, for an
algebra $A$ and an $A$-bimodule $M$, and the {\em stabilization},
a.k.a.\ {\em additivization}, a.k.a.\ the first Goodwillie
derivative of the cyclic nerve functor $|N^{cy}(F(A),F(M))|$, where
$F(A)$ is the category of finitely generated free $A$-modules, and
$F(M)$ is an $F(A)$-bimodule corresponding to $M$. If $A=k$, then it
is easy to write a functorial multiplicative map
$\phi:|N^{cy}(F(k),F(M))| \to C(M)$ to the $p$-th cyclic power $C(M)
= (M^{\otimes p})^{\Z/p\Z}$ of a vector space $M$, and then one
shows that the stabilization of $C(M)$ is the truncation $\tau^{\geq
  0}\vH^\hdot(\Z/p\Z,M^{\otimes p})$ of the Tate cohomology
algebra. Stabilizing $\phi$ gives our map $\phi_\idot$, and it
remains to show that up to a non-zero constant,
$\phi_\idot(v)=u^{-1}$. This is a computation in homological degree
$2$, so it should not be too difficult.

However, there is a catch: it turns out that it {\em is}
difficult. In the end, \cite{K.bb} resorts to a workaround: by
plugging in some additional structure --- namely, the Connes-Tsygan
differential $B$ on $THH_\idot(k)$ --- one can further reduce the
computation to something in homological degree $0$. That one,
fortunately, is indeed quite easy. However, the Connes-Tsygan
differential in question only exists as a map of spectra, and it in
fact vanishes on individual homotopy groups. This makes the
resulting argument somewhat delicate, and requires rather strong
stabilization technology to prove that everything is compatible with
everything else.

\medskip

The goal of this note is to suggest an alternative approach to
proving that the map $\phi_\idot$ above indeed sends the B\"okstedt
periodicity generator $v$ to a non-trivial multiple of $u^{-1}$. We
use a slightly different stabilization construction of $THH$:
instead of cyclic nerves, we stabilize the so-called {\em $K$-theory
  with coefficients}, as suggested and done in \cite{DM}. This is
much simpler technically since cyclic nerves typically have higher
homotopy groups, while $K$-theory with coefficients in our setting
essentially reduces to $K_0$ (``essentially'' means ``after
localization at $p$'', and this is harmless since $THH(k)$ is
$p$-local, see Lemma~\ref{p.le}). We then show that $\phi$ extends
to $K$-theory with coefficients, and exhibit an explicit degree-$2$
class $v$ in the stabilization of $K_0(k,M)$ whose image $\phi_2(v)$
is $u^{-1}$. This proves, yet again, B\"okstedt periodicity (and
shows that our explicit class is a multiplicative generator). The
construction is made possible by a couple of additional
simplifications, discussed below in Section~\ref{con.sec}. The
actual element in $K$-theory that we construct is quite simple, it
is a module of dimension $2$ (and we need to exclude the case
$p=2$). Curiously, the construction depends on an additional
parameter $\lambda \in k \setminus \{0,1\}$, and we were {\em not}
able to pinpoint a value for the parameter that makes things work;
all we prove is that such a value exists.

\medskip

The note is split into three parts. The first one is this
introduction. Then Section~\ref{gen.sec} contains all the necessary
generalities and preliminaries, and Section~\ref{con.sec} is devoted
to the actual construction. All the material in
Section~\ref{gen.sec} is standard, but for lack of convenient
references, we follow the exposition in \cite[Section 4]{K.bb} for
everything related to stabilization.

\section{Generalities.}\label{gen.sec}

\subsection{Stabilization.}

For any small category $\I$, denote by $\Ho(I)$ the category of
functors from $I$ to the category $\Top_+$ of pointed compactly
generated topological spaces localized with respect to pointwise
weak equivalences (where the localization can be constructed by any
of the standard methods, e.g.\ by model embeddings of \cite{DHKS},
and if one wishes, one can replace $\Top_+$ with the category of
pointed simplicial sets). If we let $\Ho=\Ho(\ppt)$ be the
localization of $\Top_+$ itself, then any $X \in \Ho(I)$ defines a
functor $\Ho(X):I \to \Ho$, $i \mapsto X(i)$. A functor $\gamma:I'
\to I$ between small categories induces a pullback functor
$\gamma^*:\Ho(I) \to \Ho(I')$, $X \mapsto X \circ \gamma$.

Say that a category $I$ is {\em half-additive} if it is pointed and
has finite coproducts, and for any half-additive small $I$, say that
$X \in \Ho(I)$ is {\em additive} if for any $i,i' \in I$, the
natural map
\begin{equation}\label{st.eq}
X(i \copr i') \to X(i) \times X(i')
\end{equation}
is a weak equivalence. The simplest half-additive category is the
category $\Gamma_+$ of finite pointed sets. In this case, $X \in
\Ho(\Gamma_+)$ is known as a {\em $\Gamma$-space}, and it is
additivity means being {\em special} in the sense of the Segal
machine \cite{seg}. If for any integer $n \geq 0$ we denote by $n_+
\in \Gamma_+$ the set with $n$ non-distinguished elements, then
$X(1_+)$ for any special $\Gamma$-space $X$ is a commutative monoid
object in $\Ho$, and $\pi_0(X(1_+))$ is an honest commutative
monoid. The $\Gamma$-space $X$ is called {\em group-like} if
$\pi_0(X(1_+))$ is not just a monoid but a group. For any small
category $I$, we say that $X \in \Ho(\Gamma_+ \times I)$ is {\em
  additive} resp.\ {\em stable along $\Gamma_+$} if it becomes
additive resp.\ group-like after restriction to $\Gamma_+ \times i$ for
any object $i \in I$. If $I$ is half-additive, we have the
smash-product functor $m:\Gamma_+ \times I \to I$ sending $n_+
\times i$ to the coproduct of $n$ copies of $i$, $X \in \Ho(I)$ is
additive iff $m^*X$ is additive along $\Gamma_+$, and we say that
$X$ is {\em stable} if $m^*X \in \Ho(\Gamma_+ \times I)$ is stable
along $\Gamma_+$.

For any half-additive small $I$, let $\Ho^{st}(I) \subset \Ho(I)$ be
the full subcategory spanned by stable objects. Then the embedding
$\Ho^{st}(I) \to \Ho(I)$ admits a left-adjoint {\em stabilization
  functor} $\Stab:\Ho(I) \to \Ho^{st}(I)$. Again, the stabilization
can be constructed by any of the standard methods; we follow the
overview given in \cite[Section 4]{K.bb}. Explicitly, one can extend
$X:I \to \Top_+$ to a functor $X^\Delta:\Delta^o I \to \Top_+$ from
simplicial objects in $I$ by applying it pointwise and taking
geometric realization. We then have
\begin{equation}\label{stab.eq}
\Stab(X)(i) \cong \colim_n \Omega^nX^\Delta(m(\Sigma_n \times i)),
\qquad i \in I, X \in \Ho(I),
\end{equation}
where $\Sigma_n$ is the standard simplicial $n$-sphere (that is, the
finite pointed simplicial set obtained by taking the elementary
$n$-simplex $\Delta_n$ and collapsing its $(n-1)$-skeleton to the
distinguished point). In particular, it is clear from
\eqref{stab.eq} that $\Stab$ preserves finite products and finite
homotopy limits. If $I$ has a unital monoidal structure $- \otimes
-$ that is distributive --- that is, $(i_0 \copr i_1) \otimes i'
\cong (i_0 \otimes i') \copr (i_1 \otimes i')$, and similarly in the
second variable --- then $\Ho(I)$ acquires a convolution monoidal
structure, $\Ho(X):I \to \Ho$ is lax monoidal whenever $X \in
\Ho(I)$ is an algebra object, and in this case, $\Stab(X)$ is also
an algebra object in a natural way (see \cite[Subsection
  4.5]{K.bb}).

By the Segal machine, $\Ho^{st}(\Gamma_+)$ is the category of
connective spectra. In particular, it is an additive category, so
that for any integer $n$ and object $X \in \Ho^{st}(\Gamma_+)$, we
have the endomorphism $n\id:X \to X$. If $n \geq 0$, then the
induced endomorphism $n(X):X(1_+) \to X(1_+)$ of the underlying
infinite loop space $X(1_+)$ is immediately obvious from the
$\Gamma$-space structure: it is given by the composition
\begin{equation}\label{n.eq}
\begin{CD}
X(1_+) @>{\delta}>> X(1_+)^n @>{e}>> X(n_+) @>{\mu}>> X(1_+),
\end{CD}
\end{equation}
where $\delta$ is the diagonal map, $\mu$ is induced by the
projection $n_+ \to 1_+$ sending all the non-distinguished elements
in $n_+$ to the unique non-distinguished element in $1_+$, and $e$
is inverse to the $n$-fold iteration of \eqref{stab.eq}.

For any small $I$, the full subcategory $\Ho^{st}(\Gamma_+,I)
\subset \Ho(\Gamma_+ \times I)$ spanned by objects stable along
$\Gamma_+$ is a convenient model for the homotopy category of
functors from $I$ to connective spectra. The underlying functor to
spaces is $\eps_1^*X \in \Ho(I)$, where for any $n \geq 0$, we
denote by $\eps_n:I \to \Gamma_+ \times I$ the embedding onto $n_+
\times I$. For any $i \in I$, $\eps_1(X)(i) = X(1_+ \times i)$
carries a structure of a group-like special $\Gamma$-space, and the
functor $\Ho(\eps_1^*X):I \to \Ho$ lifts to a functor $\Ho^{st}(X):I
\to \Ho^{st}(\Gamma_+)$. The category $\Ho^{st}(\Gamma_+,I)$ is
additive for any small $I$ but this requires a proof; what is
immediately obvious is a functorial endomorphism $n(X):\eps_1^*X \to
\eps_1^*X$, $n \geq 0$ obtained by replacing $X(1_+)$, $X(n_+)$ in
\eqref{n.eq} by $\eps_1^*X$, $\eps_n^*X$.

If $I$ is half-additive, we have the functor $m^*:\Ho^{st}(I) \to
\Ho^{st}(\Gamma_+, I)$. It is a fully faithful embedding, with a
left-adjoint $\Stab_I:\Ho^{st}(\Gamma_+,I) \to \Ho^{st}(I)$
isomorphic to $\Stab \circ \eps_1^*$ (see e.g.\ \cite[Subsection
  5.3]{K.bb}, and specifically \cite[Corollary
  4.11]{K.bb}). Informally, one can also stabilize functors to
connective spectra, and the underlying space of a spectral
stabilization coincides with the space-level stabilization of the
underlying space. For any stable $X$ in $\Ho(I)$ and $i \in I$,
$X(i) \cong \eps_1^*m^*X(i)$ is canonically a connective spectrum,
and if $I$ is unital monoidal and $i=1$ is the unit object, then
$X(1)$ is a ring spectrum. Stabilization $\Stab_I$ commutes with
finite products and finite homotopy limits, and by \eqref{n.eq}, for
any $X \in \Ho^{st}(\Gamma_+,I)$ and integer $n \geq 0$, we have
$\Stab(n(X)) = n(\Stab(X))$.

\begin{exa}\label{DK.exa}
An additive category $I$ is certainly half-additive, so that all the
stabilization machinery applies to $\Ho(I)$. If we further assume
that $I$ is Karoubi-closed, then we have the Dold-Kan equivalence
$\Delta^o I \cong C_{\geq 0}(I)$ between simplicial objects in $I$
and chain complexes in $I$ concentrated in non-negative homological
degrees. Then $m(\Sigma_n \times i)$ in \eqref{stab.eq} coincides
with $B^n i$, the simplicial object that corresponds to $i$
considered as a complex placed in homological degree $n$. Another
simplification is that for any stable object $X \in \Ho^{st}(I)$,
the map $n(m^*X)$ of \eqref{n.eq} is simply induced by the
endomorphism $n\id:\Id \to \Id$ of the identity functor $\Id:I \to
I$. Indeed, since $I$ is additive, the $n$-fold coproduct $m(n_+
\times i)$ of an object $i \in I$ with itself is canonically
identified with the $n$-fold product $i^n$, and then the diagonal
map $\delta$ in \eqref{n.eq} is induced by the corresponding
functorial map $\id \to m \circ \eps_n$.
\end{exa}

\subsection{Additivization.}

For any commutative ring $k$, the stabilization construction has a
$k$-linear homological counterpart. Namely, for any small category
$I$, let $\Fun(I,k)$ be the abelian category of functors from $I$ to
the category $k\amod$ of $k$-modules, let $C_{\leq 0}(I,k)$ be the
category of complexes in $\Fun(\E,k)$ concentrated in non-negative
homological degrees, and let $\Ho(I,k)$ be the category obtained by
inverting quasiisomorphisms in $C_{\leq 0}(I,k)$, so that $\Ho(I,k)$
is the connective part of the standard $t$-structure on the derived
category of $\Fun(I,k)$. Alternatively, one can apply the Dold-Kan
equivalence and interpret $C_{\leq 0}(I,k)$ as the category of
functors from $\E$ to simplicial $k$-modules, and invert pointwise
weak equivalences. We then have the forgetful functor $\Ho(I,k) \to
\Ho(I)$ sending a simplicial $k$-module to the underlying simplicial
set.

If $I$ is additive, then again, $X \in \Ho(I,k)$ is additive iff
\eqref{st.eq} is a quasiisomorphism. For an additive $X \in
\Ho(\Gamma_+,k)$, the underlying $\Gamma$-space is automatically
group-like, so that in the $k$-linear context, stability and
additivity are one and the same. The embedding $\Ho^{st}(I,k) \to
\Ho(I,k)$ of the full subcategory spanned by stable objects admits a
left-adjoint stabilization, or additivization functor $\Stab$. This
is again compatible with monoidal structures, and can be explicitly
computed by \eqref{stab.eq}, where $\Omega^n$ now stands for the
homological shift, and $X^\Delta$ is defined by taking $X$
pointwise, then applying the Dold-Kan equivalence and totalizing the
resulting bicomplex. The forgetful functor $\Ho(I,k) \to \Ho(I)$
commutes with stabilization. There are also alternative purely
homological constructions of additivization, see \cite[Subsection
  4.6]{K.bb}.

\begin{exa}\label{C.exa}
Let $k$ be a perfect field of some positive characteristic $p$, let
$I$ be the category of finite-dimensional $k$-vector spaces, and
let $C \in \Ho(I,k)$ be the {\em cyclic power functor} given by
$C(M) = H^0(\Z/p\Z,M^{\otimes p})$, $M \in I$, where $\Z/p\Z$ acts
by the longest cycle permutation. Then we have
\begin{equation}\label{C.M}
\Stab(C)(M) \cong \tau^{\leq 0}\vH^\hdot(\Z/p\Z,M^{\otimes p}),
\qquad M \in \E,
\end{equation}
where $\vH^\hdot$ stands for Tate cohomology, and $\tau^{\leq 0}$ is
the canonical truncation at $0$ (also known as the ``connective
cover''). Roughly speaking --- see \cite[Lemma 5.1]{K.bb}
for details --- to compute the right-hand side of \eqref{C.M}, one
takes a projective resolution $P_\idot$ of the trivial
$k[\Z/p\Z]$-module $k$, lets $P'_\idot$ be the cone of the
augmentation map $P_\idot \to k$, and takes $H^0(\Z/p\Z,P'_\idot
\otimes M^{\otimes p})$ (up to a quasiisomorphism, this does not
depend on the choice of $P_\idot$). Then one checks that stabilizing
$H^0(\Z/p\Z,P_\idot \otimes M^{\otimes p})$ gives $0$, and
$\tau^{\leq 0}\vH^\hdot(\Z/p\Z,M^{\otimes p})$ is stable. The cyclic
power functor $C$ is lax monoidal, so that its stabilization
$\Stab(C)$ is also lax monoidal, and $\Stab(C)(k)$ carries a natural
multiplication; this is the same multiplication as the usual
multiplication in Tate cohomology. If $p \geq 3$, then $\Stab(C)
\cong k[\eps,u]$ is the free graded-commutative algebra with two
generators $\eps$, $u$ of homological degrees $\deg\eps = 1$, $\deg
u = 2$.
\end{exa}

Note that if in Example~\ref{C.exa}, one computes $\Stab(C)$ by
\eqref{stab.eq}, then at $n$-th step, what we have is the complex
$H^0(\Z/p\Z,P(n)_\idot \otimes M^{\otimes p})$, where $P(n)_\idot$
is obtained by applying the Dold-Kand equivalence to the $p$-fold
pointwise tensor product $(B^n k)^{\otimes p}$, and shifting it by
$n$. We have $P(n)_0 \cong k$, $P(n)_l=0$ for $l > (p-1)n$, and all
the $k[\Z/p\Z]$-modules $P(n)_l$, $1 \leq l \leq (p-1)n$ are
projective. The only non-trivial homology group of the complex
$P(n)_\idot$ is $k$ in the bottom degree $(p-1)n$. When one takes
the colimit with respect to $n$, one gets the acyclic complex
$P'_\idot = \colim_nP(n)_\idot$ with $P'_0 \cong k$, and projective
$k[\Z/p\Z]$-modules $P_l$, $l \geq 1$.

\subsection{THH and K-theory.}\label{thh.subs}

Now fix a commutative ring $k$, a flat unital assocative $k$-algebra
$A$, and an $A$-bimodule $M$. Denote by $P(A)$ the category of
finitely generated projective left $A$-modules, and let $\langle
P(A),M \rangle$ be the category of pairs $\langle V,a \rangle$, $V
\in P(A)$, $a:V \to M \otimes_A V$ an arbitrary map. Then the
(split) exact structure on $P(A)$ induces an exact category
structure on $\langle P(A),M \rangle$, and one can consider the
corresponding $K$-theory space. Following \cite{DM}, we denote it
$K(A,M)$ and call it {\em $K$-theory of $A$ with coefficients in
  $M$}. The fundamental result of \cite{DM} then claims that we
have a natural homotopy equivalence
\begin{equation}\label{THH.K.eq}
THH(A,M) \cong \Stab(K)(A,M),
\end{equation}
where stabilization is taken with respect to the bimodule variable
$M$ (strictly speaking, this is not well-defined since the category
of $A$-bimodules is large, but it is obvious from \eqref{stab.eq}
that for any fixed $M$, it suffices to stablize over the small
category formed by finite sums $M^{\oplus n}$, $n \geq 0$). To
obtain individual Topological Hochschild Homology groups
$THH_\idot(A,M)$ one takes the homotopy groups $\pi_\idot(THH(A,M))$
of the space $THH(A,M)$. In the particular case $A=k$, $k$ a field,
$M$ is simply a $k$-vector space, and we have an obvious collection
of functors
\begin{equation}\label{prod.eq}
\begin{aligned}
\langle P(k),M \rangle \times \langle P(k),M' \rangle &\to \langle
P(k),M \otimes M' \rangle,\\
\langle V,a \rangle \times \langle V',a' \rangle &\mapsto \langle V
\otimes V',a \otimes a'\rangle
\end{aligned}
\end{equation}
that turns $K(k,M)$ into a lax monoidal functor; the induced
multiplication on $\Stab(K)(k,k) \cong THH(k)$ is the standard
multiplication on $THH(k)$. Note that if $M$ is finite-dimensional,
then $\langle P(k),M \rangle$ is the category of finite-dimensional
modules over the tensor algebra $T^\hdot M^*$ of the dual vector
space $M$. If we further assume that $k$ is a finite field, there is
the following aditional simplification.

\begin{lemma}\label{p.le}
Let $k$ be a finite field. Then the map $K(k,M) \to K_0(k,M)$
becomes an equivalence after stabilization with respect to $M$.
\end{lemma}

\proof{} Let $K_{\geq 1}(k,M)$ be the homotopy fiber of the
augmentation map. Since stabilization commutes with finite homotopy
limits, it suffices to prove that $K_{\geq 1}(k,-)$ has trivial
stabilization $X(-)=\Stab(K_{\geq 1})(k,-)$. Let $p=\cchar k$, and
consider the endomorphism $p(X):X \to X$ of \eqref{n.eq}. Then since
the category of $k$-vector spaces is additive and annihilated by
$p$, we have $p(X)=0$ by Example~\ref{DK.exa}, so that $X$ is a
retract of the homotopy fiber $X_p$ of the map $p(X):X \to
X$. Moreover, by their construction, $K$-theory spaces come from
connective spectra, so we can obtain $X$ by stabilizing the
corresponding object $Y(-) \in \Ho^{st}(\Gamma_+,k\amod)$ along
$k\amod$, and then $X_p$ is the stabilization of the homotopy fiber
$Y_p$ of the map $p(Y):Y \to Y$. Up to a shift, $Y_p(M) = K_{\geq
  1}(k,M)/p$ is the truncated $K$-theory of the category $\langle
P(k),M \rangle$ with coefficients in $\Z/p\Z$. This is a Noetherian
Artinian abelian category, so by devissage, we have
\begin{equation}\label{sum.eq}
K_{\geq 1}(k,M)/p \cong \bigoplus_{V} K_{\geq 1}(k(V))/p,
\end{equation}
where the sum is over all the isomorphism classes of irreducible
objects $V \in \langle p(k),M \rangle$, and $k(V)$ is the
endomorphism field of the object $V$. But $k(V)$ is a finite field
of characteristic $p$, so $K_{\geq 1}(k(V))/p=0$ by \cite{Q}.
\endproof

Lemma~\ref{p.le} shows that for a finite field, computing $THH(k,M)$
by \eqref{THH.K.eq} reduces to stabilizing $K_0(k,M)$, a functor
from $k$-vector spaces to discrete abelian groups. The values of
this functor are rather large --- as in \eqref{sum.eq}, $K_0(k,M)$
for a finite-dimensional vector space $M$ is the free abelian group
generated by all isomorphism classes of irreducible
finite-dimensional modules over $T^\hdot M^*$. Already for $\dim M =
2$ this is a countable but huge set with no discernible
structure. However, we can at least construct maps from $K_0(k,M)$
to other functors. What we need is a map
\begin{equation}\label{phi.eq}
\phi:K_0(k,M) \to C(M), \qquad \langle V,a \rangle \mapsto \Tr(a^p),
\end{equation}
where $C(M)$ is the cyclic power functor of Example~\ref{C.exa},
$a^p:V \to V \otimes M^{\otimes p}$ is the $p$-fold iteration of the
map $a$, and its the trace $\Tr(a^p) \in M^{\otimes p}$ lies in
$C(M)$ by virtue of the cyclic symmetry property of traces ($\Tr ab
= \Tr ba$). The map \eqref{phi.eq} is obviously functorial;
moreover, if we equip its source with the lax monoidal structure
induced by \eqref{prod.eq}, then it is a lax monoidal map with
respect to the standard lax monoidal structure on $C(M)$ (checking
this amounts to observing that $(a \otimes b)^p = a^p \otimes
b^p$). Therefore we can apply the stabilization functor $\Stab$ and
obtain a multiplicative map
$$
\phi_\idot = \pi_\idot(\Stab(\phi)):THH_\idot(k) \to \Stab(C)(k) \cong
k[\eps,u],
$$
where the right-hand side is as in Example~\ref{C.exa}. What we will
do is construct an element $v \in THH_2(k)$ such that $\phi_2(v) =
u$.

\section{The construction.}\label{con.sec}

\subsection{The object.}

Fix a finite field $k$ of characteristic $p = \cchar k \geq 3$. As
in Subsection~\ref{thh.subs}, let $P(k)$ be the category of
finite-dimensional $k$-vector spaces, and for any $M \in P(k)$, let
$\langle P(k),M \rangle$ be the exact category of pairs $\langle V,a
\rangle$, $V \in P(k)$, $a:V \to V \otimes_k M$. For any pair
$\langle V,a \rangle \in \langle P(k),M\rangle$, let $[V,a] \in
K_0(k,M)$ be its class. Note that if $M=k$ is one-dimensional, the
same devissage as in \eqref{sum.eq} provides an identification
\begin{equation}\label{sum.K.eq}
K_0(k,k) \cong \Z[k] \cong \bigoplus_{\lambda \in k}\Z \cdot
[k,\lambda].
\end{equation}
Consider the endofunctor $\langle P(k),M \rangle \to \langle
P(k),M\rangle$ sending $\langle V,a \rangle$ to $\langle V,-a
\rangle$, and let $\iota:K_0(k,M) \to K_0(k,M)$ be the corresponding
involution. Denote by $L(k,M) = K_0(k,M)/K_0(k,M)^{\iota}$ the
quotient of the $K_0$-group $K_0(k,M)$ by its subgroup of
$\iota$-invariants. Note that since $p$ is odd, $\tr (-a)^{\otimes
  p} = - \tr a^{\otimes p}$, so that the map \eqref{phi.eq} factors
as
\begin{equation}\label{K.io}
\begin{CD}
K_0(k,M) @>{q}>> L(k,M) @>{\psi}>> C(M),
\end{CD}
\end{equation}
where $q$ is the quotient map. Since as we saw in the proof of
Lemma~\ref{p.le}, the stabilization $\Stab(K_0)_\idot(k,M) \cong
THH_\idot(k,M)$ is $p$-local, and $p$ is odd, we have a
decomposition
\begin{equation}\label{thh.io}
THH_\idot(k,k) \cong THH_\idot(k,k)^\iota \oplus \Stab(L)_\idot(k,k),
\end{equation}
and to find an element $v \in THH_2(k,k)$ with $\phi(v) \neq 0$, it
suffices to find an element $v \in \Stab(L)_2(k,k)$ with $\psi(v)
\neq 0$ (in fact, it is easy to see from \eqref{n.eq} that $\iota$
acts on $THH_\idot(k,k)$ by $-\id$, so that $THH_\idot(k,k)^\iota$ in
\eqref{thh.io} vanishes, but we will not need it).

In order to find such an element $v$, we can compute
$\Stab(L)_\idot(k,k)$ by the colimit \eqref{stab.eq} in the form
described in Example~\ref{DK.exa}, and in fact, the desired element
$v$ appears already at the first step of the colimit.

To construct it, note that $B(k)([n]) \cong k^{\oplus n}$, $[n] \in
\Delta$, and we have the face maps $\6^i_n:B(k)([n]) \to
B(k)([n-1])$, $0 \leq i \leq n$. To simplify notation, let $L([n]) =
L(k,B(k)([n]))$, with the face maps $\6^i_n:L([n]) \to
L([n-1])$. Consider the projective line $\Pp^1(k)$ over $k$, with
the line bundle $\calo(2)$, and identify $H^0(\Pp^1(k),\calo(2))
\cong \ssl_2(k)$, the Lie algebra of vector fields on
$\Pp^1(k)$. For any point $p \in \Pp^1(k)$, evaluation at $p$ gives
a map $\ev_p:\ssl_2(k) \to k$ whose kernel $\bfr_p \subset
\ssl_2(k)$ is a Borel subalgebra. Identify $\ssl_2(k) \cong
\ssl_2(k)^*$ by the Killing form, and note that dually, we have maps
$\ev_p^*:\ssl_2(k) \to \bfr_p^*$. Let $V=k^{\oplus 2}$ be the
$2$-dimensional $k$-vector space, and let
\begin{equation}\label{al.eq}
\alpha:\ssl_2(k) \otimes_k V \to V
\end{equation}
be the standard $\ssl_2(k)$-action on $V$.

\begin{lemma}\label{borel.lm}
For any element $\lambda \in k$, $\lambda \neq 0,1$, there exists a
unique isomorphism $B(k)([2]) \cong \ssl_2(k)$ such that
$\6^i_2:B(k)([2]) \to B(k)([1]) \cong k^{\oplus 2}$ for $i=0,1,2,3$
factors through the evaluation $\ev_p^*:\ssl_2(k) \to \bfr_p$ at
$p=0,1,\infty,\lambda$. Moreover, if we let $a_\lambda:V \to V
\otimes_k B(k)([2]) \cong V \otimes_k \ssl_2(k)$ be dual to
\eqref{al.eq} under this unique isomorphism, and let $v_\lambda =
q([V,a_\lambda]) \in L([2])$ be the image of the class $[V,a] \in
K_0(k,B(k)([2]))$ under the quotient map \eqref{K.io}, then we have
$\6^i_2(v_\lambda) = 0$ for any $i=0,1,2,3$.
\end{lemma}

\proof{} The first claim is obvious: quadruples of one-dimensional
subspaces $k \subset k^{\oplus 3}$ in generic position form a single
$GL_3(k)$-orbit, and quadruples $\{\Ker \6^i_n\}$ and $\{\Ker
\ev_p^*\}$ are both in generic position (for any $\lambda \neq
0,1$). For the second claim, note that $\6^i_2([V,a])$, $i =
0,1,2,3$ is then the class $[V,b_p]$, $p = 0,1,\infty,\lambda$,
where $b_p:V \to V \otimes_k B(k)([1]) \cong V \otimes_k \bfr_p^*$ is
dual to the action $\bfr_p \otimes V \to V$ induced by
\eqref{al.eq}. Since $\bfr_p \subset \ssl_2(k)$ is a Borel subalgebra,
$V$ as a representation of $\bfr_p$ is an extension of $1$-dimensional
representations that both come by restriction via the quotient map
$u:\bfr_p \to k \cong \bfr_p/\ufr_p$, where $\ufr_p \subset \bfr_p$ is the
nilpotent radical. By devissage, we then have
$$
[V,b_p] = u^*[k,\zeta_0] + u^*[k,\zeta_1]
$$
for some elements $\zeta_0,\zeta_1 \in k$, where $[k,-]$ is as in
\eqref{sum.K.eq}. Since elements in $\bfr_p \subset \ssl_2(k)
\subset \End(V)$ are traceless, we moreover have $\zeta_1 =
-\zeta_0$, and then $[V,b_p]$ is invariant under $\iota$, so that
$\6^i_n(v) = q([V,b_p])=0$.
\endproof

\subsection{The trace computation.}
Let $\gfr=\ssl_2(k)$ with the standard basis given by the matricies
\begin{equation*}
  e = \begin{pmatrix}
    0 & 1 \\
    0 & 0 \\
  \end{pmatrix}, \qquad
  h = \begin{pmatrix}
    1 & 0 \\
    0 & -1 \\
  \end{pmatrix}, \qquad
  f = \begin{pmatrix}
    0 & 0 \\
    1 & 0 \\
  \end{pmatrix},
\end{equation*}
and let $x, y, z$ denote the basis in $\gfr^*$ dual to $e, h, f$.
We would like to explicitly write down the element
$t=\phi([V, a])$ (recall that $\phi$ was defined in \eqref{phi.eq}).
The space $C(\gfr^*)$ has a basis indexed by necklaces of length
$p$ over the alphabet ${x, y, z}$. We will use noncommutative
monomials in order to represent the corresponding basis vector.
In order to write $t$ down in the latter basis, one needs to compute
the traces of all $p$-element products of $e$, $h$, and $f$.
Consider such a product. We may assume that $e$ or $f$ occurs in 
it at least once ($p$ is odd, so $\tr h^p = \tr h = 0$).
Since $eh=e$ and $fh=-f$, dropping all the occurrences of $h$
will only change the sign. Thus, we may assume that the factors
are equal to either $e$ or $f$. Since $e^2=f^2=0$, both $e$ and $f$
are traceless, and
\begin{equation*}
  ef=e_{11} = \begin{pmatrix}
    1 & 0 \\
    0 & 0 \\
  \end{pmatrix},
\end{equation*}
we conclude that the trace of the original product is zero
unless it was of the form
\begin{equation*}
  h\ldots heh\ldots hfh\ldots heh\ldots hfh\ldots h
\end{equation*}
(in other words, in the corresponding necklace $e$ appears
at least once, every $e$ is followed by a number of $h$'s and
an $f$, and every$f$ is followed by a number of $h$'s and an $e$),
and the trace of the latter equals $1$ up to a sign.

We conclude that
\begin{equation}\label{t.eq}
  t = \sum \pm xy\ldots yzy\ldots yxy\ldots yzy\ldots y,
\end{equation}
where the explicit signs will be computed later.

\subsection{The cohomology class.}
In order to do the necessary computations, we will us the
explicit form of stabilization via the Dold-Kan equivalence;
namely, its first step.
Let $X_\bullet$ denote the simplicial group obtained from
the complex
\begin{equation*}
  k[-1] = \cdots\to 0\to k \to 0.
\end{equation*}
Denote the standard bases in $X_2$, $X_3$, and $X_4$ by $X_2=\langle w_1, w_2\rangle$,
$X_3=\langle v_1, v_2, v_3 \rangle$, and $X_4=\langle u_1, u_2, u_3, u_4 \rangle$ respectively.
The face maps $d_i: X_3\to X_2$ are given by the matrices
\begin{equation*}
  d_0 = \begin{pmatrix}
    0 & 1 & 0 \\
    0 & 0 & 1 \\
  \end{pmatrix},\
  d_1 = \begin{pmatrix}
    1 & 1 & 0 \\
    0 & 0 & 1 \\
  \end{pmatrix},\
  d_2 = \begin{pmatrix}
    1 & 0 & 0 \\
    0 & 1 & 1 \\
  \end{pmatrix},\
  d_3 = \begin{pmatrix}
    1 & 0 & 0 \\
    0 & 1 & 0 \\
  \end{pmatrix},
\end{equation*}
while the face maps $d_j:X_4\to X_3$ are given by the matrices
\begin{eqnarray*}
  d_0 = \begin{pmatrix}
    0 & 1 & 0 & 0 \\
    0 & 0 & 1 & 0 \\
    0 & 0 & 0 & 1 \\
  \end{pmatrix},\quad
  d_1 = \begin{pmatrix}
    1 & 1 & 0 & 0 \\
    0 & 0 & 1 & 0 \\
    0 & 0 & 0 & 1 \\
  \end{pmatrix},\quad
  d_2 = \begin{pmatrix}
    1 & 0 & 0 & 0 \\
    0 & 1 & 1 & 0 \\
    0 & 0 & 0 & 1 \\
  \end{pmatrix},\\
  d_3 = \begin{pmatrix}
    1 & 0 & 0 & 0 \\
    0 & 1 & 0 & 0 \\
    0 & 0 & 1 & 1 \\
  \end{pmatrix},\quad
  d_4 = \begin{pmatrix}
    1 & 0 & 0 & 0 \\
    0 & 1 & 0 & 0 \\
    0 & 0 & 1 & 0 \\
  \end{pmatrix}.
\end{eqnarray*}

We would like to find an isomorphism $\tau:\gfr^*\stackrel{\sim}{\to} X_3$
such that the face maps $d_0$, $d_1$, $d_2$, $d_3$ become dual to inclusions
of Borel subalgebras (see Lemma~\ref{borel.lm}). In other words, the inverse
images of $\Ker d_i$, $i=0,\ldots,3$, should lie on the quadric $4ef+h^2$.
A simple computation shows that
\begin{equation*}
  \begin{pmatrix} v_1 \\ v_2 \\ v_3 \\ \end{pmatrix}
  =
  \begin{pmatrix}
    \lambda(\lambda-1) & 0 & 0 \\
    \lambda^2 & 2\lambda & -\lambda \\
    0 & 0 & -(\lambda-1) \\
  \end{pmatrix}
  \begin{pmatrix} x \\ y \\ z \\ \end{pmatrix}
\end{equation*}
is a family of such morphisms indexed by $\lambda\in
\mathbb{A}^1(k)\setminus\{0, 1\}$.
Its inverse
(up to a scalar) is given by
\begin{equation}\label{alpha.eq}
  \tau_\lambda (x) = -2v_1,\quad
  \tau_\lambda (y) = v_2 + \lambda(v_1 - v_2 + v_3),\quad
  \tau_\lambda (z) = 2\lambda v_3.
\end{equation}
Extend $\tau_\lambda$ to $\lambda\in \{0,1\}$ by the same formulas,
and condsider the family of traces $t_\lambda = C(\tau_\lambda)(t)$.

\begin{prop}\label{trace.prop}
  For some $\lambda\in \mathbb{A}^1(k)\setminus\{0, 1\}$ the element
  $t_\lambda$ represents a non-trivial homology class in the normalized
  complex associated to the simplicial group $C(X_\bullet)$.
\end{prop}
The proof of this proposition will take the rest of the section.

Recall that by Lemma~\ref{borel.lm} $t_\lambda$ lies in the normalized
complex. Thus, it is enough to show that $t_\lambda$ is not in the image
of the map
\begin{equation*}
  d_0 : \bigcap_{i=1,\ldots,4} \Ker d_i \to C(X_3).
\end{equation*}

\begin{lemma}\label{t01.lm}
  One has $t_0=t_1=0$.
\end{lemma}

\proof{}
  Recall that $t$ is a linear combination of cyclic monomials of the form
  $xy\ldots yzy\ldots yxy\ldots yzy\ldots y$, and each of these monomials
  contains at least one $z$ (see~\eqref{t.eq} and the discussion preceeding
  it). Since $\tau_0(z)=0$, we immediately conclude that $t_0=0$.

  Let us now show that $t_1=C(\tau_1)(t)=0$. Recall that
  $\tau_1(x)=-2v_1$, $\tau_1(y)=v_1 + v_3$, $\tau_1(z)=2v_3$.
  It is easy to see that it is enough to prove the following.
  Let $q>1$ be an odd integer.
  Consider a length $q$ binary word of the form
  $w=0^{a_1}1^{b_1}0^{a_2}1^{b_2}\ldots 0^{a_k}1^{b_k}$
  such that $a_i, b_i > 0$ for $i=1,\ldots, k$.
  Let $P$ denote the set of sequences
  $S_w = (s_1, s_2, \ldots, s_{2l})$, $1\leq s_1 < s_2 < \cdots < s_{2l}\leq q$,
  $l > 0$,
  such that $w_{s_i}\neq w_{s_{i+1}}$ for all $i=1,\ldots, 2l-1$.
  Given a pair $(w, s)$, $s\in S_w$, let
  $p(w, s) = g_1g_2\ldots g_q$, where
  \begin{equation*}
    g_i = \begin{cases}
      h, & \text{if } i \text{ does not appear in } s, \\
      e, & \text{if } i \text{ appears in } s \text{ and } w_i=0, \\
      f, & \text{if } i \text{ appears in } s \text{ and } w_i=1,
    \end{cases}
  \end{equation*}
  and let $t(w, s) = \tr p(w,s)$. We claim that $\sum_{s\in S_w} t(w, s) = 0$.

  We prove the latter by induction on $l$. First, let $l=1$.
  Without loss of generality assume that $b_1$ is even.
  Then 
  \[S_w = \{ (s_1, s_2) \mid 1\leq s_1 \leq a_1, a_1+1 \leq s_2 \leq q \}.\]
  The involution $\iota$, which sends $(s_1, s_2)$ to $(s_1, s_2+1)$ when
  $s_2-a_1$ is odd, and to $s_2-1$ when $s_2-a_1$ is even, acts freely on $S_w$.
  Moreover, $t(w, s) = -t(w, \iota s)$ since $p(w,s)$ and $p(w,\iota s)$ differ
  by exchanging adjacent $h$ and $f$ in the product.
  
  Now assume that $l>1$. Without loss of generality assume that $b_1$ is even.
  Then $S_w = S'_w\sqcup S''_w$, where $S''_w$ denotes the set of those
  sequences containing an element from $\{a_1+1, a_1+2,\ldots, a_1+b_1\}$.
  Remark that such an element must be unique, and an argument similar to
  the one used for the base case shows that $\sum_{s\in S''_w} t(w, s) = 0$.
  It remains to observe that $S'_w$ is naturally in bijection with
  $S_{w'}$, where $w'=0^{a_1+a_2}1^{b_2}0^{a_3}1^{b_3}\ldots 0^{a_k}1^{b_k}$,
  and traces are preseved under this bijection. The claim now follows
  from the inductive hypothesis.
\endproof

Let us now look at the formal derivative $t'_\lambda(0)$.
Since, $\tau_\lambda(z)=2\lambda v_3$ and every monomial appearing
in $t$ contains at least one $z$, we conclude that the monomials
contributing to $t'_\lambda(0)$ are of the form $x y^a z y^b$, where
$a,b\geq 0$, $a+b=p-2$. Remark that the coefficient in $t$ in front
of $x y^a z y^{a}$ equals $\tr (e h^a f h^b) = (-1)^b$. Since
$\tau_\lambda(z)$ is linear homogeneous in $\lambda$,
$\tau_\lambda(y) = v_2 + \lambda (v_1-v_2+v_3)$, and
$\tau_\lambda(x)=-2v_1$, we conclude that
\begin{equation*}
  t'_\lambda(0) = -4\sum_{a=0}^{p-2}(-1)^{a}v_1v_2^{p-2-a}v_3v_2^{a}.
\end{equation*}

Before we prove the following lemma, let us introduce some
convenient notation. Let $N_k$ denote the set of length $p$
necklaces in the alphabet $\{1, 2, \ldots, a\}$. Given an element
$w=w_1w_2\ldots w_p\in N_4$, denote by $u_w$ the basis vector
in $C(X_4)$ corresponding to $u_{w_1}u_{w_2}\cdots u_{w_p}$.
Given a necklace $w\in N_3$, the basis vector $v_w$ in $C(X_3)$
is defined similarly. We abuse the notation and define
$d_2: \{1, 2, 3, 4\} \to \{1, 2, 3\}$ by the rules $d_2(1)=1$,
$d_2(2)=d_2(3)=2$, $d_2(4)=3$. Then for any $w\in N_4$ one
has $d_2(u_w)=v_{d_2(w)}$. The maps $d_1$ and $d_3$ are defined
similarly.

We will need the following subsets of $N_k$. Let $N^{1,1}_k$
denote the set of those necklaces which contain exactly one
bead equal to $1$ and exactly one bead equal to $k$.
Let $N^{+, 1}_k$ denote the set of those necklaces which
contain more than one
bead equal to $1$ and exactly one bead equal to $k$.
The set $N^{1,+}_k$ is defined similarly.
Finally, the set $N^{+,+}_k$ consists of those necklaces which
contain more that one bead equal to $1$ and more than one bead
equal to $k$.
One has the following identities:
\begin{equation}\label{filt.eq}
\begin{split}    
  d_2(N^{?, ?}_4) \subseteq N^{?,?}_3,\quad
  & d_1(N^{1,1}_4) \subseteq N^{1,1}_3\sqcup N^{+,1}_3,\quad
  d_1(N^{+,1}_4) \subseteq N^{+,1}_3, \\
  & d_3(N^{1,1}_4) \subseteq N^{1,1}_3\sqcup N^{1,+}_3,\quad
  d_3(N^{1,+}_4) \subseteq N^{1,+}_3.
\end{split}
\end{equation}
Finally, denote by $\bar{N}_k=N^{1,1}_k\sqcup N^{1,+}_k
\sqcup N^{+,1}_k \sqcup N^{+,+}_k$ (these are all the necklaces
containing at least one bead equal to $1$ and at least one bead
equal to $k$).

\begin{lemma}\label{tprime.lm}
  The element $t'(0)$ is not in the image of the map
  \[
    d_0 : \bigcap_{i=1,\ldots,4} \Ker d_i \to C(X_3).
  \]
\end{lemma}

\proof{} For convenience, we will work with $v = -t'_\lambda(0)/4$
instead. Precisely,
\begin{equation*}
  v = v_1v_2^{p-2}v_3 - v_1v_2^{p-3}v_3v_2 + \cdots +
  v_1v_2v_3v_2^{p-3} - v_1v_3v_2^{p-2}.
\end{equation*}

Consider the element $u'\in C(X_4)$ given by
\begin{equation*}
  \begin{split}
    u' & = u_2u_3^{p-2}u_4 - u_2u_3^{p-3}u_4u_3 + \cdots +
    u_2u_3u_4u_3^{p-3} - u_2u_4u_3^{p-2} - \\
    & - u_1u_3^{p-2}u_4 + u_1u_3^{p-3}u_4u_3 - \cdots -
    u_1u_3u_4u_3^{p-3} + u_1u_4u_3^{p-2}.
  \end{split}
\end{equation*}
It satisfies the equations $d_0(u')=v$, $d_1(u')=0$,
$d_2(u')=-v$, $d_3(u')=0$, $d_4(u')=0$. Thus, it is enough to show
that there is no solution to the system of equations
\begin{equation}
d_0(u) = d_1(u) = d_3(u) = d_4(u) = 0,\quad d_2(u)=v.
\end{equation}

Consider the subspace $\Ker d_0$. Since $d_0(u_0) = 0$ and
$d_0(u_i) = v_{i-1}$ for $i=1,2,3$, we see that $\Ker d_0$
is spanned by the basis vectors $u_w$ for those $w\in N_4$
which contain at least one bead equal to $1$. Similarly,
$\Ker d_4$ is spanned by the basis vectors corresponding
to necklaces containing at least one bead equal to $4$.
We conclude that $\bar{U} = \Ker d_0 \cap \Ker d_4$ is
spanned by the basis vectors $u_w$, where $w$ runs
through the set $\bar{N}_4$.


Recall that $U=\langle u_w \mid u_w \in \bar{N}_4\rangle$,
and put $V=\langle v_w \mid v_w \in \bar{N}_3\rangle$.
It follows from the identities~\eqref{filt.eq} that
$d_i(U)\subseteq V$ for $i=1,2,3$. Since $v\in V$, it is
enough to show that there exists no $u\in U$ such that
$d_2(u)=v, d_1(u)=d_3(u)=0$. In order to do that, we are
going to construct an element $\nu\in V^*$ such that
$\nu(v)\neq 0$, while $d^*_2(\nu)$ equals zero when
restricted to $\Ker d_2 \cap \Ker d_3$.

Let $\langle u^\vee_w \mid w\in \bar{N}_4\rangle$
and $\langle v^\vee_w \mid w\in \bar{N}_3\rangle$ denote the bases
in $U^*$ and $V^*$ dual to
$\langle u_w \mid w\in \bar{N}_4\rangle$ and
$\langle v_w \mid w\in \bar{N}_3\rangle$ respectively.
Considet the map $p:X_4\to [1]$ given by $p(1)=p(2)=0$,
$p(3)=p(4)=1$. Let $w=w_1w_2\ldots w_p$ be a word in
the alphabet $\{1, 2, 3, 4\}$ with at least one occurrence
of $1$ and at least one occurrence of $4$. Let $0\leq i < p$
be such that $\sigma^i p(w)$ is lexicographically minimal,
where $\sigma$ denotes the cyclic permutation.
If
$1\leq i_1 < i_2 < \ldots < i_a\leq p$ are the indicies
of the letters $w_i$ in $w$ equal to $1$, put
\begin{equation*}
  \gamma_1(w) = \frac{i_1+i_2+\cdots+i_a}{a} - i \in k.
\end{equation*}
Similarly, if $1\leq j_1 < j_2 < \ldots < j_b\leq p$ are
the indicies of the letters $w_j$ in $w$ equal to $4$,
put 
\begin{equation*}
  \gamma_4(w) = \frac{j_1+j_2+\cdots+j_b}{b} - i \in k.
\end{equation*}
Finally, put $\Delta(w) = \gamma_1(w)-\gamma_4(w)$. Remark that
both $s_1$ and $s_4$ are invariant under the cyclic
action. Thus, $s_1$, $s_4$, and $\Delta$ descend to
maps from $\bar{N}_4$ to $k$.

Consider the element $\mu\in U^*$ given by
\begin{equation*}
  \mu = \sum_{w\in \bar{N}_4} \Delta(w)u^\vee_w.
\end{equation*}
We claim that $\mu(\Ker d_1\cap \Ker d_3)=0$.
Indeed, pick sections $s_i:\bar{N}_3\to\bar{N}_4$
of the maps $d_i$, $i=1,2,3$. Then
\begin{equation*}
  \mu = \sum_{w\in \bar{N}_3}
  \left(\gamma_1(s_1(w))d_1^* - \gamma_4(s_3(w))d_3^*\right) v^\vee_w,
\end{equation*}
and the left hand side trivially vanishes on $\Ker d_1\cap \Ker d_3$.
Meanwhile, since $\delta(w)$ depends only on the positions
the beads equal to $1$ and $4$, we see that $\mu=d_2^*\nu$,
where
\begin{equation*}
  \nu = \sum_{w\in \bar{N}_3}
  \Delta(s_2(w)) v^\vee_w.
\end{equation*}

Finally, let us compute $\nu(v)$. For $w\in N^{1,1}_4$
the number $\delta$ measures the clockwise distance from the
bead equal to $1$ to the bead equal to $4$. Thus
\[
  \nu(v) = 1-2+3-4+\ldots + (p-2) - (p-1) = \frac{1-p}{2}\neq 0.
\]
\endproof

The proof of Proposition~\ref{trace.prop} now trivially follows from
Lemmas~\ref{t01.lm} and~\ref{tprime.lm} and the fact that the parametric
curve $t_\lambda$ is of degree at most $p-1$.


{\small\noindent
Affiliations (for both authors):
\begin{enumerate}
\renewcommand{\labelenumi}{\arabic{enumi}.}
\item Steklov Mathematics Institute (main affiliation).
\item National Research University Higher School of Economics.
\end{enumerate}}

{\small\noindent
{\em E-mail addresses\/}: {\tt kaledin@mi-ras.ru,} {\tt avfonarev@mi-ras.ru.}
}

\end{document}